\documentclass[journal]{IEEEtran}
\usepackage{amsmath,amsfonts,amsthm}
\usepackage{algorithmic}
\usepackage{algorithm}
\usepackage{array}
\usepackage[caption=false,font=normalsize,labelfont=sf,textfont=sf]{subfig}
\usepackage{textcomp}
\usepackage{stfloats}
\usepackage{url}
\usepackage{verbatim}
\usepackage{graphicx}
\usepackage{cite}
\usepackage[american, americanvoltages]{circuitikz}
\usepackage{bm}
\usepackage[hidelinks]{hyperref}
\theoremstyle{remark}
\newtheorem*{remark}{Remark}

\begin{document}

\title{Optimal Sensitivity of the general\\ Wheatstone Bridge}

\author{Michael Fischer
  \thanks{Michael Fischer is an independent researcher based in Fulda, Germany
  (e-mail: \href{mailto:michael.fischer.research@gmail.com}
  {michael.fischer.research@gmail.com}).}}

\maketitle

\begin{abstract}
Optimizing the sensitivity of the unbalance voltage in Wheatstone bridges with respect to changes in the bridge parameters remains a fundamental objective in circuit design and instrumentation. When finite source and detector resistances are taken into account, determining the optimal bridge configuration becomes increasingly complex, and a closed-form analytical representation of the optimal solution has not yet been established. This paper derives an analytical solution for the optimal configuration of a Wheatstone bridge with finite source and detector resistances. Furthermore, the proposed optimal solution is compared with the conventional equal-arm configuration.
\end{abstract}


\section{Introduction}
\IEEEPARstart{T}{he} Wheatstone bridge is a well-established circuit topology for precise resistance measurements and resistive sensor interfacing \cite{Frank1959}. Originally developed by Samuel Hunter Christie in 1833 and popularized by Charles Wheatstone in 1843 \cite{Ekelof2001}, the Wheatstone bridge remains a core component in modern precision instrumentation systems to map relative resistance variations—such as those from strain gauges or piezoresistive transducers—directly into a differential output voltage \cite{Pallas2001}.

Maximizing the sensitivity of the unbalance voltage with respect to resistance variations is thus essential for improving measurement precision in such applications. Consequently, finding the bridge configuration yielding optimal sensitivity represents a fundamental objective in circuit design. While the optimal conditions for the ideal Wheatstone bridge with zero source resistance and infinite detector resistance are well understood, the problem becomes non-trivial when accounting for finite source and detector resistances. In fact, when source and detector resistances are of the same order of magnitude as the bridge arms, the general Wheatstone bridge must be considered \cite{Weiss1969}.

Back in 1917, Laws already dedicated a section entitled "Sensitiveness Attainable with the Wheatstone Bridge" in his textbook
\cite{Laws1917}, primarily focusing on the thermal limitations of the bridge components. In 1959, Frank utilized the compensation theorem to derive an approximate analytical representation of the null sensitivity as a function of the remaining bridge components and derived the maximum null sensitivity for an equal-arm setup \cite{Frank1959}. In 1969, Weiss introduced the dimensionless sensitivity function and built on the results from \cite{Frank1959} to derive analytical results for the null sensitivity of the general Wheatstone bridge. To this end, Weiss introduces the term "universal curves" in order to elegantly express the null sensitivity as a function of the involved components. However, Weiss does not provide an analytical representation of the optimal Wheatstone bridge configuration and suggests the equal-arm design as a practical choice. Since then, numerous further contributions regarding the optimization of the Wheatstone bridge sensitivity have been made \cite{Maisel1977, Takagishi1978, Takagishi1980, Torrents2025}. However, to the best of the authors' knowledge, the overall optimal configuration has not yet been analytically derived. In the recent past, the problem has been approached computationally by numerically approximating optimal configurations for given finite source and detector resistances \cite{Zhou2025}.

In the present paper, an analytical solution for the optimal Wheatstone bridge configuration with finite source and detector resistances is derived. Furthermore, an analytical representation of the corresponding null sensitivity is provided and benchmarked against the null sensitivity of the equal-arm setup. Finally, an explicit functional relationship between the change in the unbalance voltage for a given change in the measured resistance value is established.

\section{\texorpdfstring{Wheatstone bridge setup and \\ problem definition}{Wheatstone bridge setup and problem definition}}
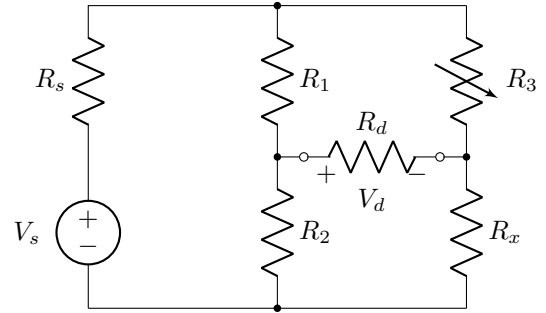
\begin{figure}[!t]
  \centering
  \begin{circuitikz}[every node/.style={transform shape=false}]
    
    % Quelle und Innenwiderstand R_s
    \draw (0,0)
      to[vsource, v=\(V_s\), invert] (0,2)
      to[R, l=\(R_s\)] (0,4)
      -- (2.5,4);
    
    % Obere Brückenverbindung
    \draw (2.5,4) -- (5,4);
    
    % Linker Brückenzweig
    \draw (2.5,4)
      to[R, l=\(R_1\)] (2.5,2)
      to[R, l=\(R_2\)] (2.5,0);
          
    % Rechter Brückenzweig
    \draw (5,4)
      to[variable resistor, l=\(R_3\)] (5,2)
      to[R, l=\(R_x\)] (5,0);
          
    % Detektor-Zweig (Messbrücke)
    
    % Linke Anschlussleitung und Klemme
    \draw (2.5,2) -- (2.81,2);
    \draw (2.85,2) node[ocirc, inner sep=0pt] {};
    \draw (2.89,2) -- (3.2,2);

    % Detektorwiderstand
    \draw (3.2,2)
      to[R, l^=\(R_d\), v_>=\(V_d\)] (4.3,2);

    % Rechte Anschlussleitung und Klemme
    \draw (4.3,2) -- (4.61,2);
    \draw (4.65,2) node[ocirc, inner sep=0pt] {};
    \draw (4.69,2) -- (5,2);
          
    % Untere Brückenverbindung
    \draw (0,0) -- (5,0);

    % Knotenpunkte
    \path (5,2) node[circle, fill, inner sep=1pt] {};
    \path (2.5,4) node[circle, fill, inner sep=1pt] {};
    \path (2.5,2) node[circle, fill, inner sep=1pt] {};
    \path (2.5,0) node[circle, fill, inner sep=1pt] {};

  \end{circuitikz}

  \caption{Circuit diagram of the general Wheatstone bridge with fixed \(R_1,R_2\), adjustable \(R_3\) and unknown resistance \(R_x\)}
  \label{fig:Wheatstone_circuit}
\end{figure}

We consider the general Wheatstone bridge with source and detector resistances \(0<R_s<\infty\) and \(0<R_d<\infty\), respectively. Without loss of generality, \(R_x\) denotes the unknown resistor to be measured, \(R_1\) and \(R_2\) are fixed, and \(R_3\) is adjustable, cf. Fig.~\ref{fig:Wheatstone_circuit}. This is a standard setup in the literature and coincides with case \(1\) analyzed by Weiss in \cite{Weiss1969}. It can be shown, that the alternative configuration with a variable \(R_2\) and a fixed \(R_3\) is equivalent to interchanging source and detector \cite{Weiss1969}. Many other equivalent bridge configurations are already investigated in \cite{Frank1959} and are therefore not in the scope of this paper. 

The bridge is said to be balanced when the voltage across \(R_d\) vanishes, i.e., \(V_{d}=0,\) which is equivalent to the null condition \(R_1R_x=R_2R_3.\) In this setup, for given resistances \(R_1\) and \(R_2,\) the resistance \(R_3\) is adjusted such that \(V_{d}=0.\) Then, the unknown resistance \(R_x\) can be obtained by \(R_{x} = R_{3}\frac{R_{2}}{R_{1}}.\) Furthermore, with \(k=\frac{R_{1}}{R_{2}},\) the adjustable resistance \(R_3\) is implicitly given by \(R_3=kR_x\) in the balanced state.

As in \cite{Weiss1969}, the dimensionless sensitivity function \(S_{R_{i}}\) is defined by the fractional change in resistance \(R_i\) of the normalized, unbalanced voltage \(V_{d},\) i.e.,
\begin{equation}\label{eq: Sensi}
  S_{R_{i}} = \frac{\frac{\partial V_{d}}{V_s}}{\frac{\partial R_{i}}{R_{i}}} = \frac{R_{i}}{V_s}\frac{\partial V_{d}}{\partial R_{i}}=\frac{R_{i}R_d}{V_s}\frac{\partial I_{d}}{\partial R_{i}}
\end{equation}
for \(i=1,2,3,x.\) Furthermore, let \(\overline{S_{R_{i}}}\) denote the corresponding null sensitivity, i.e., \(S_{R_{i}}\) with satisfied null condition \(R_1R_x=R_2R_3.\) Without loss of generality, the remaining analysis is done with respect to the sensitivity \(S_{R_{x}}.\) It will also be shown that all results can easily be transferred to the alternative sensitivity definitions, such as \(S_{R_{2}}\) or \(S_{R_{3}}.\)

Since \(R_3\) is adjusted such that the system is balanced, the overall setup is uniquely determined by \(R_1,R_2,R_x\) as well as the source and detector resistances \(R_{s}\) and \(R_d,\) respectively. Since \(R_1=k R_2,\) the overall setup and the null sensitivity \(\overline{S_{R_{x}}}\) are thus uniquely determined by \(k,\ R_{2},\ R_{x},\) as well as the source and detector resistances \(R_{s}\) and \(R_d,\) respectively.

Thus, for given source and detector resistances \(R_{s}\) and \(R_d,\) the question arises, which choice of \(R_1\) and \(R_2,\) or equivalently which values of \(k\) and \(R_{2}\) maximize the null sensitivity \(\overline{S_{R_{x}}}\) for a given \(R_{x}.\) This results in the following two dimensional optimization problem for the null sensitivity:
\begin{equation}\label{eq: optimization problem}
  \max_{\left( k,R_{2} \right) \in \mathbb{R}_{+}^{2}}\left| \overline{S_{R_{x}}}\left( k,R_{2},R_{x} \right) \right| = :\left| \overline{S_{R_{x}}^{opt}}\left( R_{x} \right) \right|
\end{equation}

Furthermore, this raises the additional question, which value of resistance \(R_{x}\) determines the global optimum of the null sensitivity \(\overline{S_{R_{x}}^{opt}}\left( R_{x} \right).\)

To the best of the authors' knowledge, the optimal pair of resistances \(R_1^{opt},R_{2}^{opt}\) or equivalently \(k^{opt},R_{2}^{opt}\) for given source and detector resistances \(R_{s}\) and \(R_d,\) and a given value of \(R_{x},\) has not yet been derived in the literature.

\section{Derivation of the null sensitivity}
The circuit diagram of the general Wheatstone bridge shown in Fig.~\ref{fig:Wheatstone_circuit} contains \(4\) nodes and \(6\) branches. Thus, \(4-1=3\) linearly independent node equations and \(6-4+1=3\) independent loop equations are required for a unique characterization of the system, respectively, cf. \cite{Feldmann1991}. The resulting system of equations is given by
\begin{align*}
  I_{s}R_{s} + I_{1}R_{1} + \left( I_{1} - I_{d} \right)R_{2} - V_s &= 0 \\
  \left( I_{s} - I_{1} \right)R_{3} - I_{d}R_{d} - I_{1}R_{1} &= 0 \\
  \left( I_{s} - I_{1} + I_{d} \right)R_{x} - \left( I_{1} - I_{d} \right)R_{2} + I_{d}R_{d} &= 0,
\end{align*}
which can equivalently be represented via \(\bm{R}\cdot\bm{I}=\bm{V}\) with
\begin{equation}\label{eq: R}
  \bm{R}=
  \left(
    \begin{array}{ccc}
      R_{s} & R_{1} + R_{2} & - R_{2} \\
      - R_{3} & R_{1} + R_{3} & R_{d} \\
      - R_{x} & R_{2} + R_{x} & - \left( R_{2} + R_{d} + R_{x} \right)
    \end{array}
  \right),
\end{equation}
\(\bm{I}=\left(
    \begin{array}{c}
      I_{s},~I_{1},~I_{d}
    \end{array}
  \right)^T\)  
and
\(\bm{V}=\left(
    \begin{array}{c}
      V_s,~0,~0
    \end{array}
  \right)^T.\)
Note that a similar representation has also been chosen by \cite{Zhou2025}.

For the remainder of this paper and any given quantity \(X\), let \(\overline{X}\) denote \(X\) at the balanced state, i.e., with satisfied null condition \(R_1R_x=R_2R_3,\) or equivalently with \(R_1=kR_2,R_3=kR_x.\) Thus, \(\overline{\bm{R}}\) is given by
\begin{equation}\label{eq: Rbalanced}
  \overline{\bm{R}}=
  \left(
    \begin{array}{ccc}
      R_{s} & (1 + k)R_{2} & - R_{2} \\
      - {kR}_{x} & k(R_{2} + R_{x}) & R_{d} \\
      - R_{x} & R_{2} + R_{x} & - \left( R_{2} + R_{d} + R_{x} \right)
    \end{array}
  \right).
\end{equation}

For the remainder of the analysis, the determinant \(\det(\overline{\bm{R}})\) will play a major role. For an elegant calculation of \(\det(\overline{\bm{R}}),\) the second row of \(\overline{\bm{R}}\) is replaced by the second row minus \(k\) times the third row of \(\overline{\bm{R}}\), which results in
\begin{equation*}
  \widetilde{\bm{R}}=
  \left(
    \begin{array}{ccc}
      R_{s} & (1 + k)R_{2} & - R_{2} \\
      0 & 0 & R_{d}+k\left( R_{2} + R_{d} + R_{x} \right) \\
      - R_{x} & R_{2} + R_{x} & - \left( R_{2} + R_{d} + R_{x} \right)
    \end{array}
  \right)
\end{equation*}
with \(\det(\overline{\bm{R}})=\det(\widetilde{\bm{R}}).\) Thus, \(\det(\overline{\bm{R}})\) can easily be evaluated via
\begin{equation}
\label{detBalancedR}
\begin{split}
  \det(\overline{\bm{R}}) = &-\left( R_{d} + k\left( R_{2} + R_{d} + R_{x} \right) \right) \\
                       &\quad\cdot\left( \left( R_{s} + (1 + k)R_{x} \right)R_{2} + R_{s}R_{x} \right).
\end{split}
\end{equation}

As a result, the overall network behavior, represented by \(\det(\overline{\bm{R}}),\) is factorized into two separate products, each having a distinct physical interpretation. The second product is entirely independent of the detector resistance \(R_d.\) A short calculation shows that this product describes the isolated network behavior at balanced state, where instead of including \(R_d\) the detector branch is completely omitted, i.e., unloaded and treated as an open circuit. On the other hand, the first product includes \(R_{d}\) alongside the core bridge parameters, thereby representing the coupling of the detector to the overall bridge network. Although no current flows through \(R_{d}\) at exact balance, this term dictates the bridge sensitivity near the balanced state. It determines how a minor perturbation in \(R_x\) translates into a change in \(V_{d}\) or \(I_{d}\) at the detector.

For an elegant derivation of \(\overline{S_{R_{x}}},\) it will be convenient to express \(\det(\overline{\bm{R}})\) as a polynomial in both, \(k\) and \(R_2:\)
\begin{equation}
\label{detBalancedRPolyk}
\begin{split}
  \det(\overline{\bm{R}})= &- R_{2}R_{x}\left( R_{2} + R_{d} + R_{x} \right)\cdot k^{2}\\
                      &- \biggl( \left( R_{2} + R_{d} + R_{x} \right)\left( \left( R_{s} + R_{x} \right)R_{2} + R_{s}R_{x} \right)\\
                      &\quad+ R_{d}R_{2}R_{x}\biggr)\cdot k\\
                      & - R_{d}\left( \left( R_{s} + R_{x} \right)R_{2} + R_{s}R_{x} \right)
\end{split}
\end{equation}
\begin{equation}
\label{detBalancedRPolyR2}
\begin{split}
  \det(\overline{\bm{R}})= &- k\left( R_{s} + (1 + k)R_{x} \right)\cdot R_{2}^{2}\\
                      &- \biggl(\left( R_{d} + k\left( R_{d} + R_{x} \right) \right)\left( R_{s} + (1 + k)R_{x} \right)\\
                      &\quad+kR_sR_x\biggr)\cdot R_{2}\\
                      & - \left( R_{d} + k\left( R_{d}+{R}_{x} \right) \right)R_{s}R_{x}
\end{split}
\end{equation}

By utilizing Cramer's rule, the current \(I_d\) can be evaluated via
\begin{align}\label{eq: Id_Cramer}
  I_d=\frac{\det(\bm{R_V})}{\det(\bm{R})},
\end{align}
where \(\bm{R_V}\) denotes the Cramer matrix formed by replacing the last column of \(\bm{R}\) by
\(\bm{V}=\left(
  \begin{array}{c}
    V_s,~0,~0
  \end{array}
\right)^T.\)
Thus, with the representation of \(\bm{R}\) given in (\ref{eq: R}), one obtains
\begin{equation}\label{eq: CramerDet}
  \begin{split}
  \det(\bm{\bm{R_V}})&=V_s\left( -R_{3}(R_{2} + R_{x}) + (R_{1} + R_{3})R_{x}\right) \\
           &= V_s({R_{1}R_{x} - R_{2}R}_{3}),
  \end{split}
\end{equation}
which again yields \(\det\left(\overline{\bm{R_{V}}}\right) = 0.\) Furthermore, it follows 
\begin{align*}
  \frac{\partial}{\partial R_{x}}\det\left(\bm{R_{V}}\right) = V_sR_{1}
\end{align*}
and thus
\begin{equation}\label{eq:partialDetRUBalanced}
  \overline{\frac{\partial}{\partial R_{x}}\det\left( \bm{R_{V}} \right)} = kV_{s}R_{2}
\end{equation}
for the balanced state. According to (\ref{eq: Sensi}), the sensitivity \(S_{R_{x}}\) is given by \(S_{R_{x}}=\frac{R_{d}R_x}{V_s}\frac{\partial I_{d}}{\partial R_{x}}.\) Therefore, by utilizing (\ref{eq: Id_Cramer}), the sensitivity \(S_{R_{x}}\) can be expressed by
\begin{align*}
  S_{R_{x}}&=\frac{R_{d}R_x}{V_s}\frac{\partial}{\partial R_{x}}\left( \frac{\det(\bm{R_{V}})}{\det(\bm{R})}\right)\\
           &=\frac{R_{d}R_{x}}{V_s}\\
           &\quad\cdot\frac{\det(\bm{R})\frac{\partial}{\partial R_{x}}\det\left( \bm{R_{V}} \right) - \det\left( \bm{R_{V}} \right)\frac{\partial}{\partial R_{x}}\det(\bm{R})}{{\det(\bm{R})}^{2}}.
\end{align*}
While this representation might appear unnecessarily complicated at first glance, by leveraging \(\det\left(\overline{\bm{R_{V}}}\right) = 0\) and (\ref{eq:partialDetRUBalanced}), it immediately implies the very compact representation
\begin{equation}\label{eq: SRx_compact}
  \begin{split}
    \overline{S_{R_{x}}}\left( k,R_{2},R_{x} \right) &= \frac{R_{d}R_{x}}{V_s}\frac{\overline{\frac{\partial}{\partial R_{x}}\det\left( \bm{R_{V}} \right)}}{\det\left( \overline{\bm{R}} \right)}\\
                                                     &= \frac{kR_{2}R_{d}R_{x}}{\det\left( \overline{\bm{R}} \right)}
   \end{split}
\end{equation}
for the null sensitivity. By inserting the evaluation of \(\det\left(\overline{\bm{R}}\right)\) given in (\ref{detBalancedR}), one obtains the explicit representation
\begin{align*}
  &\overline{S_{R_{x}}}\left( k,R_{2},R_{x} \right)= \\
  &\frac{- kR_{2}R_{d}R_{x}}{\left( R_{d} + k\left( R_{2} + R_{d} + R_{x} \right) \right)\left( \left( R_{s} + (1 + k)R_{x} \right)R_{2} + R_{s}R_{x} \right)}
\end{align*}
for the null sensitivity for given source and detector resistances \(R_s\) and \(R_d\), respectively.

Since the network determinant \(\det\left( \overline{\bm{R}} \right)\) appears in the denominator of \(\overline{S_{R_{x}}}\left( k,R_{2},R_{x} \right),\) the null sensitivity is inherently governed by both the detector-to-bridge coupling and the open-circuit network behavior, which are described by the first and second factors of \(\det\left( \overline{\bm{R}} \right),\) respectively.
\begin{remark}
  Completely analogous to the above derivation of \(\overline{S_{R_{x}}}\left( k,R_{2},R_{x} \right),\) it follows
\begin{align*}
  \overline{S_{R_{3}}}\left( k,R_{2},R_{x} \right) &= \frac{R_{d}R_{3}}{V_s}\frac{\overline{\frac{\partial}{\partial R_{3}}\det\left( \bm{R_{V}} \right)}}{\det\left( \overline{\bm{R}} \right)} \\
                                                   &= \frac{{- R}_{2}R_{d}R_{3}}{\det\left( \overline{\bm{R}} \right)} \\
                                                   &= \frac{{- kR}_{2}R_{d}R_{x}}{\det\left( \overline{\bm{R}} \right)} = - \overline{S_{R_{x}}}\left( k,R_{2},R_{x} \right)
\end{align*}
for the alternative null sensitivity \(\overline{S_{R_{3}}}\left( k,R_{2},R_{x}\right)\) with respect to the variable resistance \(R_3,\) which is frequently used in the literature, cf. e.g., case \(1\) in \cite{Weiss1969}. Since \( \overline{S_{R_{3}}}=- \overline{S_{R_{x}}},\) all results derived in the following, can immediately be transferred to this alternative definition of sensitivity.

Similarly, sensitivity results for \(\overline{S_{R_{2}}}\) with variable \(R_2\) and constant resistance values \(R_1\) and \(R_3\), corresponding to case 2 in \cite{Weiss1969}, can be derived from \(\overline{S_{R_{x}}}.\)
\end{remark}

\section{\texorpdfstring{Derivation of the optimal \\ Wheatstone bridge configuration}{Derivation of the optimal Wheatstone bridge configuration}}
In this section, exact solutions \(k^{opt}\) and \(R_2^{opt}\) of the optimization problem (\ref{eq: optimization problem}) are being derived, resulting in an optimized null sensitivity \(\overline{S_{R_{x}}^{opt}}\left( R_{x} \right)=\overline{S_{R_{x}}}\left( k^{opt},R_{2}^{opt},R_{x} \right)\) for a given value \(R_x\) and given source and detector resistances \(R_s\) and \(R_d,\) respectively. To this end, equations for the partial derivatives \(\frac{\partial}{\partial k}\overline{S_{R_{x}}}\) and \(\frac{\partial}{\partial R_2}\overline{S_{R_{x}}}\) are being derived and set to zero.

With the help of (\ref{eq: SRx_compact}), the partial derivatives with respect to \(k\) and \(R_2\) satisfy
\begin{align*}
  \frac{\partial}{\partial k}\overline{S_{R_{x}}}\left( k,R_{2},R_{x} \right)&= R_{2}R_{d}R_{x}\frac{\det\left( \overline{\bm{R}} \right) - k\frac{\partial}{\partial k}{\det\left( \overline{\bm{R}} \right)}}{{\det\left( \overline{\bm{R}} \right)}^{2}}
\end{align*}
and
\begin{align*}
  \frac{\partial}{\partial R_2}\overline{S_{R_{x}}}\left( k,R_{2},R_{x} \right)&=kR_{d}R_{x}\frac{\det\left( \overline{\bm{R}} \right) - R_{2}\frac{\partial}{\partial R_{2}}{\det\left( \overline{\bm{R}} \right)}}{{\det\left( \overline{\bm{R}} \right)}^{2}},
\end{align*}
respectively. Now, the representations of \(\det\left( \overline{\bm{R}} \right)\) as polynomials in \(k\) and \(R_2\), cf. (\ref{detBalancedRPolyk}) and (\ref{detBalancedRPolyR2}), can be leveraged effectively. Since \(p(x)-x\frac{d}{d x}p(x)=-ax^2+c\) for any polynomial \(p(x)=ax^2+bx+c,\) (\ref{detBalancedRPolyk}) and (\ref{detBalancedRPolyR2}) yield
\begin{align*}
  &\frac{\partial}{\partial k}\overline{S_{R_{x}}}\left( k,R_{2},R_{x} \right)= R_{2}R_{d}R_{x}\\
  &\cdot
  \frac{R_{2}R_{x}\left( R_{2} + R_{d} + R_{x} \right)k^{2}-R_{d}\left( \left( R_{s} + R_{x} \right)R_{2} + R_{s}R_{x} \right)}{{\det\left( \overline{\bm{R}} \right)}^{2}}
\end{align*}
and
\begin{align*}
  &\frac{\partial}{\partial R_2}\overline{S_{R_{x}}}\left( k,R_{2},R_{x} \right)=\\
  &kR_{d}R_{x}\frac{k\left( R_{s} + (1 + k)R_{x} \right)R_{2}^{2} - \left( R_{d} + k\left( R_{d}{+ R}_{x} \right) \right)R_{s}R_{x}}{{\det\left( \overline{\bm{R}} \right)}^{2}}
\end{align*}
for the partial derivatives. Therefore, \(\frac{\partial}{\partial k}\overline{S_{R_{x}}}=0\) if and only if
\begin{align}\label{eq: stationaryK}
  k^2\left( R_{2},R_{x} \right) = \frac{R_{d}\left( \left( R_{s} + R_{x} \right)R_{2} + R_{s}R_{x} \right)}{R_{2}R_{x}\left( R_{2} + R_{d} + R_{x} \right)}.
\end{align}
Analogously, \(\frac{\partial}{\partial R_2}\overline{S_{R_{x}}}=0\) if and only if
\begin{align}\label{eq: stationary R2}
  R_{2}^2\left( k,R_{x} \right) =\frac{R_{s}R_{x}\left( R_{d} + k\left( R_{d}{+ R}_{x} \right) \right)}{k\left( R_{s} + (1 + k)R_{x} \right)}.
\end{align}

With the help of (\ref{eq: SRx_compact}) or the subsequent explicit representation of \(\overline{S_{R_{x}}}\left( k,R_{2},R_{x} \right),\) it is easily observed that \(\overline{S_{R_{x}}}(k,R_2,R_x)\) equals zero for \(k=0\) or \(R_2=0\) and vanishes asymptotically as \(\|(k,R_2)\|\to\infty.\) Since \(\overline{S_{R_{x}}}(k,R_2,R_x)<0\) for all \((k,R_2)>0,\) a standard compactness argument yields that \(\overline{S_{R_{x}}}\) must attain its global minimum \((k^{opt},R_2^{opt})\) in the interior of \(\mathbb{R}^2_+.\) Since \(\overline{S_{R_{x}}}\) is continuously differentiable with respect to \((k,R_2),\) all minima must be stationary points, characterized by \(\nabla_{(k,R_2)}\overline{S_{R_{x}}}=0\) or equivalently by the nonlinear system of equations given by (\ref{eq: stationaryK}) and (\ref{eq: stationary R2}).

A standard approach for solving the system of equations given by (\ref{eq: stationaryK}) and (\ref{eq: stationary R2}) would be to substitute one variable into the equation given by the other one. However, since both equations involve linear and quadratic dependencies on the other variable, the resulting equation would become highly nonlinear, involving polynomials of higher order, algebraic fractions and radicals in the target variable, making an analytical closed-form solution impossible. Nonetheless, the system of equations conceals a particular algebraic structure that enables the derivation of a unique analytical solution by leveraging the following trick.

Instead of substituting \(R_{2}\left( k,R_{x} \right)\) given by (\ref{eq: stationary R2}) into all linear and quadratic \(R_2\)-dependencies in 
\begin{align}\label{eq: stationaryKlinquad}
  k^2\left( R_{2},R_{x} \right)=\frac{R_{d}\left( \left( R_{s} + R_{x} \right)R_{2} + R_{s}R_{x} \right)}{R_{2}^{2}R_{x} + R_{2}R_{x}\left( R_{d} + R_{x} \right)},
\end{align}
only the quadratic dependency \(R_2^2\) is substituted. This results in the scalar equation
\begin{align*}
  k^2=\frac{R_{d}\left( \left( R_{s} + R_{x} \right)R_{2} + R_{s}R_{x} \right)}{\frac{R_{s}R_{x}^2\left( R_{d} + k\left( R_{d}{+ R}_{x} \right) \right)}{k\left( R_{s} + (1 + k)R_{x} \right)} + R_{2}R_{x}\left( R_{d} + R_{x} \right)}
\end{align*}
with two variables \(k\) and \(R_2.\) Expanding this equation yields
\begin{align*}
  0 &= k^{2}\frac{R_{s}R_{x}^{2}\left( R_{d} + k\left( R_{d} + R_{x} \right) \right)}{k\left( R_{s} + (1 + k)R_{x} \right)} \\
    &\quad+ k^{2}\left( R_{d} + R_{x} \right)R_{2}R_{x} - R_{d}\left( \left( R_{s} + R_{x} \right)R_{2} + R_{s}R_{x} \right),
\end{align*}
which can be further expanded to the equivalent expression
\begin{align*}
  0 &= kR_{s}R_{x}^{2}\left( R_{d} + k\left( R_{d} + R_{x} \right) \right)\\
    &\quad+\left( R_{s} + (1 + k)R_{x} \right)\Bigl( k^{2}\left( R_{d} + R_{x} \right)R_{2}R_{x}\\
                                               &\qquad- R_{d}\left( \left( R_{s} + R_{x} \right)R_{2} + R_{s}R_{x} \right) \Bigr).
\end{align*}

Now, a careful rearrangement of the involved terms yields the equivalent relation
\begin{align*}
  0 &= kR_{s}R_{x}^{2}\left( R_{d} + k\left( R_{d} + R_{x} \right) \right)\\
    &\quad+ \left( R_{s} + (1 + k)R_{x} \right)R_2\Bigl( k^{2}\left( R_{d} + R_{x} \right)R_{x}\\
    &\qquad - R_{d}\left( R_{s} + R_{x} \right)\Bigr) -\left( R_{s} + (1 + k)R_{x} \right) R_{s}{R_{d}R}_{x}\\
    &=\left( R_{s} + (1 + k)R_{x} \right)R_2\Bigl( k^{2}\left( R_{d} + R_{x} \right)R_{x}\\
    &\qquad - R_{d}\left( R_{s} + R_{x} \right)\Bigr)\\
    &\quad+R_{s}R_{x}\Bigl(kR_{x}\left( R_{d} + k\left( R_{d} + R_{x} \right) \right)\\
    &\quad\qquad-\left( R_{s} + (1 + k)R_{x} \right)R_{d}\Bigr)\\
    &=\left( R_{s} + (1 + k)R_{x} \right)R_2\Bigl( k^{2}\left( R_{d} + R_{x} \right)R_{x}\\
    &\qquad - R_{d}\left( R_{s} + R_{x} \right)\Bigr)\\
    &\quad+R_{s}R_{x}\Bigl( k^{2}\left( R_{d} + R_{x} \right)R_{x} - R_{d}\left( R_{s} + R_{x} \right) \Bigr)\\
    &=\Bigl( \left( R_{s} + (1 + k)R_{x} \right)R_{2} + R_{s}R_{x} \Bigr)\\
    &\quad\qquad\cdot\Bigl( k^{2}\left( R_{d} + R_{x} \right)R_{x} - R_{d}\left( R_{s} + R_{x} \right) \Bigr).
\end{align*}
Since \(\left( R_{s} + (1 + k)R_{x} \right)R_{2} + R_{s}R_{x}>0\) for arbitrary \(R_2>0,\) the previous relation yields the optimality condition
\begin{align}\label{eq: optimalityConditionK}
  k^{2}\left( R_{d} + R_{x} \right)R_{x} - R_{d}\left( R_{s} + R_{x} \right)=0
\end{align}
for \(k,\) and therefore the unique optimal solution
\begin{equation}\label{eq: kOpt}
  k^{opt}\left( R_{x} \right) = \sqrt{\frac{R_{d}(R_{s} + R_{x})}{R_{x}\left( R_{d} + R_{x} \right)}}
\end{equation}
for given source and detector resistances, \(R_s\) and \(R_d,\) respectively.

Notably, the key to eliminating the additional degree of freedom \(R_2\) from the optimality condition, ultimately yielding a unique optimal solution for \(k\), is to carve out the term \(\left( R_{s} + (1 + k)R_{x} \right)R_{2} + R_{s}R_{x},\) representing the open-circuit network behavior, as a common factor in the optimality condition.

For the derivation of the global optimum \(R_2^{opt}(R_x),\) it is convenient to expand (\ref{eq: stationaryKlinquad}) via
\begin{equation}\label{eq: optimalConditionR2}
\begin{split}
  0 &= k^{2}R_{2}^{2}R_{x} + k^{2}R_{2}R_{x}\left( R_{d} + R_{x} \right) \\
    &\quad - R_{d}\left( \left( R_{s} + R_{x} \right)R_{2} + R_{s}R_{x} \right) \\
    &= k^{2}R_{x}R_{2}^{2} - R_{s}R_{d}R_{x}\\
    &\quad + \left( k^{2}R_{x}\left( R_{d} + R_{x} \right) - R_{d}\left( R_{s} + R_{x} \right) \right)R_{2}.
\end{split}
\end{equation}
This expression also contains the factor \(k^{2}R_{x}\left( R_{d} + R_{x} \right) - R_{d}\left( R_{s} + R_{x} \right)\) already known from (\ref{eq: optimalityConditionK}). Since \(k^{opt}\left( R_{x} \right)\) was chosen exactly such that \(k^{2}R_{x}\left( R_{d} + R_{x} \right) - R_{d}\left( R_{s} + R_{x} \right)=0,\) the optimal solution \(R_2^{opt}(R_x)\) can also be analytically derived from (\ref{eq: optimalConditionR2}) and is given by 
\begin{equation}\label{eq: R2Opt}
  R_{2}^{opt}(R_{x}) = \sqrt{\frac{{R_{s}R_{d}R}_{x}}{{{(k}^{opt}(R_{x}))}^{2}R_{x}}} = \sqrt{\frac{R_{s}R_{x}\left( R_{d} + R_{x} \right)}{R_{s} + R_{x}}}
\end{equation}
for given source and detector resistances, \(R_s\) and \(R_d,\) respectively.

Note that  \(k^{opt}\left( R_{x} \right)\) and \( R_{2}^{opt}(R_{x})\) from (\ref{eq: kOpt}) and (\ref{eq: R2Opt}) already uniquely determine the entire optimal Wheatstone bridge configuration. In fact, \(R_{1}^{opt}(R_{x})\) can directly be derived via
\begin{equation}\label{eq: R1Opt}
  R_{1}^{opt}\left( R_{x} \right) = k^{opt}\left( R_{x} \right)R_{2}^{opt}\left( R_{x} \right) = \sqrt{R_{s}R_{d}}.
\end{equation}
The remaining optimal resistance \(R_{3}^{opt},\) which—unlike the fixed choices of \(R_{1}^{opt}\) and \(R_{2}^{opt}\)—is implicitly determined by the balancing process, can analogously be derived by
\begin{equation*}
  R_{3}^{opt}\left( R_{x} \right) = k^{opt}\left( R_{x} \right)R_{x} = \sqrt{\frac{R_{d}R_{x}\left( R_{s} + R_{x} \right)}{R_{d} + R_{x}}}
\end{equation*}
for given source and detector resistances, \(R_s\) and \(R_d,\) respectively.

It is worth mentioning that the optimal resistance \(R_1^{opt}\) is completely independent of \(R_x\) and equals the geometric mean \(\sqrt{R_{s}R_{d}}\) of the source and detector resistances. In contrast to that, \(R_2^{opt}\) and \(R_3^{opt}\) are dependent on \(R_x.\) Moreover, \(R_2^{opt}\) and \(R_3^{opt}\) exhibit a symmetry in such a way that \(R_2^{opt}(R_x,R_s,R_d)=R_3^{opt}(R_x,R_d,R_s)\) for arbitrary source and detector resistances. Finally, one observes that the necessary null condition
\begin{align*}
  R_{2}^{opt}R_{3}^{opt} &= \sqrt{\frac{R_{s}R_{x}\left( R_{d} + R_{x} \right)}{R_{s} + R_{x}}}\sqrt{\frac{R_{d}R_{x}\left( R_{s} + R_{x} \right)}{R_{d} + R_{x}}} \\
                         &= \sqrt{R_{s}R_{d}}R_{x} \\
                         &= R_{1}^{opt}R_{x}
\end{align*}
is satisfied.

\section{Evaluation of the optimal null sensitivity}
According to (\ref{eq: SRx_compact}), an evaluation of the optimal null sensitivity \(\overline{S_{R_{x}}}\left( k^{opt},R_{2}^{opt},R_{x} \right)\) requires the evaluation of \(\det(\overline{\bm{R}}^{opt}),\) i.e., (\ref{detBalancedR}) evaluated at \(\left(k^{opt},R_2^{opt}\right).\) The resulting expression is explicitly calculated in (\ref{eq: detOpt}).

\begin{figure*}[!t]
  \normalsize
\begin{equation}\label{eq: detOpt}
  \begin{split}
    -\det(\overline{\bm{R}}^{opt})&=\left( R_{d} + k^{opt}\left( R_{2}^{opt} + R_{d} + R_{x} \right) \right)\Bigl(\left( R_{s} + (1 + k^{opt})R_{x} \right)R_{2}^{opt} + R_{s}R_{x} \Bigr)\\
  &=\left( R_{d} + \sqrt{R_{s}R_{d}} + k^{opt}\left( R_{d} + R_{x} \right) \right)\left( \sqrt{R_{s}R_{d}}R_{x} + \left( R_{s} + R_{x} \right)R_{2}^{opt} + R_{s}R_{x} \right)\\
  &=\left( R_{d} + \sqrt{R_{s}R_{d}} + \sqrt{\frac{R_{d}(R_{s} + R_{x})\left( R_{d} + R_{x} \right)}{R_{x}}} \right)\Bigl( \sqrt{R_{s}R_{d}}R_{x} + \sqrt{R_{s}R_{x}\left( R_{d} + R_{x} \right)\left( R_{s} + R_{x} \right)}+ R_{s}R_{x} \Bigr)\\
                                  &=\sqrt{R_{s}R_{d}R_{x}}\left( \sqrt{R_{d}} + \sqrt{R_{s}} + \sqrt{\frac{(R_{s} + R_{x})\left( R_{d} + R_{x} \right)}{R_{x}}} \right)\Bigl( \sqrt{R_{d}R_{x}} + \sqrt{\left( R_{d} + R_{x} \right)\left( R_{s} + R_{x} \right)}+ \sqrt{R_{s}R_{x}} \Bigr)\\
    &=\sqrt{R_{s}R_{d}}\Bigl( \sqrt{R_{s}R_{x}} + \sqrt{\left( R_{s} + R_{x} \right)\left( R_{d} + R_{x} \right)}+\sqrt{R_{d}R_{x}} \Bigr)^{2}
   \end{split}                                                                                                      
\end{equation}
\hrulefill 
\vspace*{4pt} 
\end{figure*}

With the help of (\ref{eq: SRx_compact}), the null sensitivity for the optimal Wheatstone setup is therefore given by
\begin{equation}\label{eq: nullSensiOpt}
\begin{split}
  &\overline{S_{R_{x}}^{opt}}\left( R_{x} \right) = \\
  &\quad\frac{{- R}_{d}R_{x}}{\left( \sqrt{R_{s}R_{x}} + \sqrt{\left( R_{s} + R_{x} \right)\left( R_{d} + R_{x} \right)} + \sqrt{{R_{d}R}_{x}} \right)^{2}}
\end{split}
\end{equation}
for given source and detector resistances, \(R_s\) and \(R_d,\) respectively.

The question arises as to which choice of \(R_x\) causes the optimal sensitivity \(\overline{S_{R_{x}}^{opt}}\left( R_{x} \right)\) to attain its peak value. For the calculation of this optimal value \(R_x^{opt},\) one observes that
\begin{align*}
  \overline{S_{R_{x}}^{opt}}\left( R_{x} \right) = \frac{{- R}_{d}}{\left( \sqrt{R_{s}} + \sqrt{R_{d}} + \sqrt{R_{s} + R_{d} + R_{x} + \frac{R_{s}R_{d}}{R_{x}}} \right)^{2}}.
\end{align*}
Therefore, \(\overline{S_{R_{x}}^{opt}}\) attains its optimum at the precise \(R_x\) that minimizes \(R_{x} + \frac{R_{s}R_{d}}{R_{x}}.\)
The inequality of arithmetic and geometric means gives
\[R_{x} + \frac{R_{s}R_{d}}{R_{x}} \geq 2\sqrt{R_{x}\frac{R_{s}R_{d}}{R_{x}}} = 2\sqrt{R_{s}R_{d}},\]
where equality holds if and only if \(R_{x} = \frac{R_{s}R_{d}}{R_{x}},\) i.e., if and only if \(R_{x} = \sqrt{R_{s}R_{d}}.\) As a consequence, \(\overline{S_{R_{x}}^{opt}}\) attains its optimal sensitivity at
\begin{align*}
  R_x^{opt}=\sqrt{R_{s}R_{d}}.
\end{align*}

The corresponding value of overall optimal sensitivity can be derived via
\begin{align*}
  \overline{S_{R_{x}}^{opt}}\left( R_{x}^{opt} \right) &= \frac{{- R}_{d}}{\left( \sqrt{R_{s}} + \sqrt{R_{d}} + \sqrt{R_{s} + R_{d} + 2\sqrt{R_{s}R_{d}}} \right)^{2}}\\
                                                       &= \frac{{- R}_{d}}{\left( \sqrt{R_{s}} + \sqrt{R_{d}} + \sqrt{\left( \sqrt{R_{s}} + \sqrt{R_{d}} \right)^{2}} \right)^{2}}\\
                                                       &= \frac{{- R}_{d}}{{4\left( \sqrt{R_{s}} + \sqrt{R_{d}} \right)}^{2}}\\
                                                       &=\frac{-1}{{4\left(1+ \sqrt{\frac{R_{s}}{R_d}}\right)}^{2}}.
\end{align*}
This immediately implies the global sensitivity limit \[\left| \overline{S_{R_{x}}} \right|\leq \left| \overline{S_{R_{x}}^{opt}}\right| \leq\frac{1}{{4\left(1+ \sqrt{\frac{R_{s}}{R_d}}\right)}^{2}}< \frac{1}{4}\] and furthermore shows that the optimal value of sensitivity depends on the ratio \(\sqrt{\frac{R_{s}}{R_d}}\) of source and detector resistance.

Finally, with \(R_x^{opt}=\sqrt{R_{s}R_{d}}\) and (\ref{eq: kOpt}) it follows
\begin{align*}
  k^{opt}\left( R_{x}^{opt} \right) &= \sqrt{\frac{R_{d}(R_{s} + \sqrt{R_{s}R_{d}})}{\sqrt{R_{s}R_{d}}\left( R_{d} + \sqrt{R_{s}R_{d}} \right)}}\\
  &=\sqrt{\frac{\sqrt{R_s}R_{d}(\sqrt{R_s} + \sqrt{R_{d}})}{\sqrt{R_{s}}R_{d}\left( \sqrt{R_{d}} + \sqrt{R_{s}} \right)}}=1,
\end{align*}
which in combination with (\ref{eq: R1Opt}) therefore ultimately yields
\begin{align}\label{eq: optimalEqualArm}
  R_1^{opt}=R_2^{opt}(R_x^{opt})=R_3^{opt}(R_x^{opt})=R_x^{opt}=\sqrt{R_{s}R_d}.
\end{align}
This shows that the optimal setup converges to an equal-arm design—characterized by \(k=1\)—whenever \(R_x\) approaches \(R_x^{opt}=\sqrt{R_{s}R_d}.\)

\begin{remark}
  In (\ref{eq: Sensi}), the dimensionless sensitivity function \(S_{R_x}\) is chosen and in the course of this paper, a Wheatstone bridge configuration is derived that optimizes precisely this definition of sensitivity. Since \(\frac{\partial V_{d}}{\partial R_{x}}=S_{R_x}\cdot\frac{V_s}{R_x},\) the derived sensitivity results for \(S_{R_x}\) can easily be transferred to the alternative, non-dimensionless sensitivity function \(\frac{\partial V_{d}}{\partial R_{x}},\) which is defined with respect to absolute instead of relative changes in \(R_x.\)

It can be verified that exactly the same optimal solutions \(k^{opt}(R_x)\) and \(R_2^{opt}(R_x)\) as in (\ref{eq: kOpt}) and (\ref{eq: R2Opt}), respectively, also apply for the alternative sensitivity function \(\frac{\partial V_{d}}{\partial R_{x}}.\) However, (\ref{eq: nullSensiOpt}) immediately implies that the corresponding optimal null sensitivity value, given by \(\overline{S_{R_{x}}^{opt}}\left( R_{x} \right)\cdot\frac{V_s}{R_x},\) attains its optimum at \(R_x=0,\) resulting in a physically non-meaningful configuration.
\end{remark}

\section{Sensitivity comparison of the optimal design vs. the classical equal-arm design}\label{sec: comparison}
With the help of the general representation of the null sensitivity \(\overline{S_{R_{x}}}\left( k,R_{2},R_{x} \right)\) given in (\ref{eq: SRx_compact}), it is possible to compare different bridge configurations with the optimal setup with respect to their respective null sensitivity. In fact, let
\[P(k,R_2,R_x)=\frac{\overline{S_{R_{x}}}\left( k,R_{2},R_{x} \right)}{\overline{S_{R_{x}}^{opt}}\left( R_{x} \right)}\in[0,1]\]
denote the measure of null sensitivity performance of any given bridge setup characterized by \(k\) and \(R_2\) relative to the optimal setup.

\begin{figure}[htbp]
    \centering
    \includegraphics[scale=1]{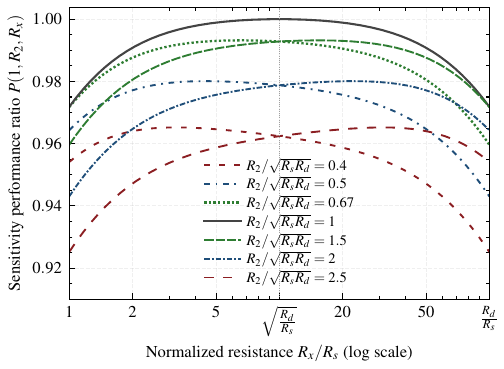}
    \caption{Sensitivity performance ratio \(P(1,R_2, R_x)\) as a function of the normalized resistance \(R_x/R_s\) for varying values of \(R_2\) and \(\frac{R_d}{R_s}=100\).}
    \label{fig:equal_arm_bridge}
\end{figure}

In the literature, the most frequently considered Wheatstone bridge configuration is the equal-arm design, which is characterized by \(k=1.\) According to (\ref{detBalancedR}), (\ref{eq: SRx_compact}) and (\ref{eq: nullSensiOpt}), the null sensitivity performance of the equal-arm design can explicitly be evaluated by
\begin{align*}
  &P(1,R_2,R_x)=\\
  &\quad\frac{R_{2}\left( \sqrt{R_{s}R_{x}} + \sqrt{\left( R_{s} + R_{x} \right)\left( R_{d} + R_{x} \right)} + \sqrt{{R_{d}R}_{x}} \right)^{2}}{\Bigl( 2R_{d} + R_{2} + R_{x} \Bigr)\Bigl( \left( R_{s} + 2R_{x} \right)R_{2} + R_{s}R_{x} \Bigr)}.
\end{align*}
For different values of \(R_2\) and \(\frac{R_d}{R_s}=100,\) the null sensitivity performance of the equal-arm design is visualized in Fig.~\ref{fig:equal_arm_bridge}.

In Fig.~\ref{fig:equal_arm_bridge}, it stands out that two equal-arm setups with
\(R_{2} = \alpha\sqrt{R_{s}R_{d}}\) and \(R_{2} = \frac{1}{\alpha}\sqrt{R_{s}R_{d}}\) for \(0 < \alpha \leq 1\)
exhibit a symmetry with respect to their null sensitivity at \(R_{x} = \sqrt{R_{s}R_{d}}.\) This symmetry is explicitly derived for arbitrary \(0 < \alpha,\beta \leq 1\) in (\ref{eq: SensiSymmetry}).

\begin{figure*}[!t]
  \normalsize
\begin{equation}\label{eq: SensiSymmetry}
  \begin{split}
    \overline{S_{R_{x}}}\left(1,\frac{1}{\alpha}\sqrt{R_{s}R_{d}},\frac{1}{\beta}\sqrt{R_{s}R_{d}} \right) &= \frac{- \frac{1}{\alpha}\sqrt{R_{s}R_{d}}R_{d}\frac{1}{\beta}\sqrt{R_{s}R_{d}}}{\left( 2R_{d} + \frac{1}{\alpha}\sqrt{R_{s}R_{d}} + \frac{1}{\beta}\sqrt{R_{s}R_{d}} \right)\left( \left( R_{s} + 2\frac{1}{\beta}\sqrt{R_{s}R_{d}} \right)\frac{1}{\alpha}\sqrt{R_{s}R_{d}} + R_{s}\frac{1}{\beta}\sqrt{R_{s}R_{d}} \right)} \\
                                                                                                           &= \frac{- \alpha\sqrt{R_{s}R_{d}}R_{d}\beta\sqrt{R_{s}R_{d}}}{\left( 2\alpha\beta R_{d} + \beta\sqrt{R_{s}R_{d}} + \alpha\sqrt{R_{s}R_{d}} \right)\left( \left( R_{s} + 2\frac{1}{\beta}\sqrt{R_{s}R_{d}} \right)\beta\sqrt{R_{s}R_{d}} + R_{s}\alpha\sqrt{R_{s}R_{d}} \right)} \\
                                                                                                           &= \frac{- \alpha\sqrt{R_{s}R_{d}}R_{d}\beta\sqrt{R_{s}R_{d}}}{\left( 2\alpha\beta R_{d} + \beta\sqrt{R_{s}R_{d}} + \alpha\sqrt{R_{s}R_{d}} \right)R_{s}\left( \beta\sqrt{R_{s}R_{d}} + 2R_{d} + \alpha\sqrt{R_{s}R_{d}} \right)} \\
                                                                                                           &= \frac{- \alpha\sqrt{R_{s}R_{d}}R_{d}\beta\sqrt{R_{s}R_{d}}}{\left( 2R_{d} + \alpha\sqrt{R_{s}R_{d}} + \beta\sqrt{R_{s}R_{d}} \right)R_{s}\left( \alpha\sqrt{R_{s}R_{d}} + 2\alpha\beta R_{d} + \beta\sqrt{R_{s}R_{d}} \right)}\\
                                                                                                           & = \frac{- \alpha\sqrt{R_{s}R_{d}}R_{d}\beta\sqrt{R_{s}R_{d}}}{\left( 2R_{d} + \alpha\sqrt{R_{s}R_{d}} + \beta\sqrt{R_{s}R_{d}} \right)\Bigl( \left( R_{s} + 2\beta\sqrt{R_{s}R_{d}} \right)\alpha\sqrt{R_{s}R_{d}} + R_{s}\beta\sqrt{R_{s}R_{d}} \Bigr)}\\
                                                                                                           & = \overline{S_{R_{x}}}\left(1,\alpha\sqrt{R_{s}R_{d}},\beta\sqrt{R_{s}R_{d}} \right)
   \end{split}                                                                                                      
\end{equation}
\hrulefill 
\vspace*{4pt} 
\end{figure*}

Furthermore, Fig.~\ref{fig:equal_arm_bridge} indicates that among all equal-arm setups, the version with \(R_2=\sqrt{R_sR_d}\) exhibits the highest null sensitivity performance. However, it reaches the performance of the optimal setup only at \(R_x=\sqrt{R_sR_d}\) with \(P(1,\sqrt{R_sR_d},\sqrt{R_sR_d})=1.\) All remaining equal-arm variants exhibit a progressively lower null sensitivity performance as \(R_2\) deviates further from \(\sqrt{R_sR_d}.\) Overall, it can be observed that, especially as \(R_x\to R_s\) or \(R_x\to R_d,\) all equal-arm setups suffer a loss in sensitivity compared with the optimal Wheatstone bridge configuration.

Finally, Fig.~\ref{fig:equal_arm_bridge} suggests that not \(R_s\) and \(R_d\), but rather their ratio \(r^2=\frac{R_d}{R_s}\) is the relevant factor affecting the null sensitivity performance. The visualizations in Fig.~\ref{fig:equal_arm_bridge} are based a ratio between detector and source resistance of \(r^2=\frac{R_d}{R_s}=100.\) Hence, the question arises as to how \(\frac{R_d}{R_s}\) affects the null sensitivity performance \(P(1,R_2, R_x).\)

To this end, let \(R_2=\alpha\sqrt{R_sR_d}\) and \(R_x=\beta\sqrt{R_sR_d}\) for \(0<\alpha,\beta\leq 1.\) Hence, with \(r=\sqrt{\frac{R_d}{R_s}}\) it holds \(\frac{R_2}{R_s}=\alpha r,\frac{R_x}{R_s}=\beta r\) and thus
\begin{align*}
  &P\left(1, R_{2},R_{x} \right) \\
  &=\frac{R_{2}\left( \sqrt{R_{s}R_{x}} + \sqrt{\left( R_{s} + R_{x} \right)\left( R_{d} + R_{x} \right)} + \sqrt{{R_{d}R}_{x}} \right)^{2}}{\Bigl( 2R_{d} + R_{2} + R_{x} \Bigr)\Bigl( \left( R_{s} + 2R_{x} \right)R_{2} + R_{s}R_{x} \Bigr)} \\
  &=\frac{\frac{R_{2}}{R_s}\left( \sqrt{\frac{R_{x}}{R_{s}}} + \sqrt{\left( 1 + \frac{R_{x}}{R_{s}} \right)\left( \frac{R_{d}}{R_{s}} + \frac{R_{x}}{R_{s}} \right)} + \sqrt{\frac{R_{d}}{R_{s}}\frac{R_{x}}{R_{s}}} \right)^{2}}{\left( 2\frac{R_{d}}{R_{s}} + \frac{R_{2}}{R_s} + \frac{R_{x}}{R_{s}} \right)\left( \left( 1 + 2\frac{R_{x}}{R_{s}} \right)\frac{R_{2}}{R_s} + \frac{R_{x}}{R_{s}} \right)}\\
  &=\frac{\alpha r\left( \sqrt{\beta r} + \sqrt{\left( 1 + \beta r \right)\left( r^2 + \beta r \right)} + \sqrt{r^2\beta r} \right)^{2}}{\left( 2r^2 + \alpha r + \beta r \right)\left( \left( 1 + 2\beta r \right)\alpha r + \beta r \right)},
\end{align*}
resulting in 
\begin{align*}
  &P\left(1, R_{2},R_{x} \right)\\
  &\quad= \frac{\alpha\left( \sqrt{\beta r} + \sqrt{(1 + \beta r)\left( r^{2} + \beta r \right)} + \sqrt{\beta r^{3}} \right)^{2}}{\left( 2r^{2} + \alpha r + \beta r \right)(\alpha + 2\alpha\beta r + \beta)}\\
                               &\underset{r \rightarrow \infty}{\sim}\frac{{\alpha\left( 2\sqrt{\beta r^{3}} \right)}^{2}}{4\alpha\beta r^{3}} = 1.
\end{align*}
This shows that for any equal-arm design defined by \(R_2=\alpha\sqrt{R_sR_d}\) and any fixed value \(R_x=\beta\sqrt{R_sR_d}\) (where \(0<\alpha,\beta\leq 1\)), the corresponding null sensitivity performance ratio \(P\left(1, R_{2},R_{x} \right)\) converges to \(1\) from below as \(r=\sqrt{\frac{R_d}{R_s}}\to\infty.\) In the context of Fig.~\ref{fig:equal_arm_bridge}, this implies that each individual curve on the logarithmic x-axis eventually approaches \(1\) from below as \(r=\sqrt{\frac{R_d}{R_s}}\rightarrow\infty.\)

In the following, the  rate of convergence of \(P\left(1, R_{2},R_{x} \right)\) to \(1\) is determined. To this end, let
\begin{align*}
  N(\alpha,\beta,r)&=\alpha\left( \sqrt{\beta r} + \sqrt{(1 + \beta r)\left( r^{2} + \beta r \right)} + \sqrt{\beta r^{3}} \right)^{2}\\
  &\quad-\left( 2r^{2} + \alpha r + \beta r \right)(\alpha + 2\alpha\beta r + \beta)
\end{align*}
and
\begin{align*}
  D(\alpha,\beta,r)=\left( 2r^{2} + \alpha r + \beta r \right)(\alpha + 2\alpha\beta r + \beta).
\end{align*}
Then, it holds true that
\begin{align*}
  P(1,R_2,R_x)-1=\frac{N(\alpha,\beta,r)}{D(\alpha,\beta,r)}.
\end{align*}

Since Taylor expansion gives
\[\sqrt{\beta + x} = \sqrt{\beta} + \frac{x}{2\sqrt{\beta}} - \frac{x^{2}}{8\beta^{\frac{3}{2}}}\mathcal{+ O}\left( x^{3} \right)\]
as \(x \rightarrow 0,\) the middle square root term in \(N(\alpha,\beta,r)\) satisfies
\begin{align*}
  \sqrt{(1 + \beta r)\left( r^{2} + \beta r \right)} &= r^{\frac{3}{2}}\sqrt{\beta + \frac{1 + \beta^{2}}{r} + \frac{\beta}{r^{2}}} \\
                                                     &= \sqrt{\beta}r^{\frac{3}{2}} + \frac{1 + \beta^{2}}{2\sqrt{\beta}}r^{\frac{1}{2}}\mathcal{+ O}\left( r^{- \frac{1}{2}} \right)
\end{align*}
as \(r\to\infty.\) Therefore, it holds true that
\begin{align*}
  &\sqrt{\beta r} + \sqrt{(1 + \beta r)\left( r^{2} + \beta r \right)} + \sqrt{\beta r^{3}} \\
  &\quad= \sqrt{\beta}r^{\frac{1}{2}} + 2\sqrt{\beta}r^{\frac{3}{2}} + \frac{1 + \beta^{2}}{2\sqrt{\beta}}r^{\frac{1}{2}}\mathcal{+ O}\left( r^{- \frac{1}{2}} \right) \\
  &\quad= 2\sqrt{\beta}r^{\frac{3}{2}} + \frac{(\beta + 1)^{2}}{2\sqrt{\beta}}r^{\frac{1}{2}}\mathcal{+ O}\left( r^{- \frac{1}{2}} \right)
\end{align*}
and thus
\begin{align*}
  &\left( \sqrt{\beta r} + \sqrt{(1 + \beta r)\left( r^{2} + \beta r \right)} + \sqrt{\beta r^{3}} \right)^{2} \\
  &\quad= 4\beta r^{3} + 2(\beta + 1)^{2}r^{2} + \mathcal{O}(r).
\end{align*}
For the entire numerator \(N(\alpha,\beta,r),\) this implies
\begin{align*}
  N(\alpha,\beta,r) &= \alpha\left( 4\beta r^{3} + 2(\beta + 1)^{2}r^{2} + \mathcal{O}(r) \right) \\
                    &\quad- \left( 2r^{2} + \alpha r + \beta r \right)(2\alpha\beta r + \alpha + \beta) \\
                    &= 4\alpha\beta r^{3} + 2\alpha(\beta + 1)^{2}r^{2} \\
                    &\quad- 4\alpha\beta r^{3} - 2(\alpha + \beta)(\alpha\beta + 1)r^{2} + \mathcal{O}(r) \\
                    &= - 2\beta(\alpha - 1)^{2}r^{2} + \mathcal{O}(r).
\end{align*}
On the other hand, the denominator \(D(\alpha,\beta,r)\) satisfies 
\begin{align*}
  D(\alpha,\beta,r)&=\left( 2r^{2} + \alpha r + \beta r \right)(\alpha + 2\alpha\beta r + \beta)\\
                   &=4\alpha\beta r^3+\mathcal{O}(r^2),
\end{align*}
which ultimately implies
\begin{align*}
  P(1,R_2,R_x)-1&=\frac{N(\alpha,\beta,r)}{D(\alpha,\beta,r)}\\
                &= \frac{- 2\beta(\alpha - 1)^{2}r^{2} + \mathcal{O}(r)}{4\alpha\beta r^{3} + \mathcal{O}\left( r^{2} \right)}\\
                &= \frac{- 2\beta(\alpha - 1)^{2}r^{- 1} + \mathcal{O}\left( r^{- 2} \right)}{4\alpha\beta + \mathcal{O}\left( r^{- 1} \right)}\\
                &= -\frac{(\alpha - 1)^{2}}{2\alpha}\frac{1}{r} + \mathcal{O}\left( r^{-2} \right).
\end{align*}

As a result, it holds true that
\begin{equation*}
{1-P(1,R_2,R_x)} = \begin{cases}
  \mathcal{O}\left( r^{-1} \right),&{\text{if}}\ {\alpha\neq1} \\ 
{\mathcal{O}\left( r^{-2} \right),}&{\text{if}}\ {\alpha=1} 
\end{cases}
\end{equation*}
for \(r=\sqrt{\frac{R_d}{R_s}}\) with \(R_2=\alpha\sqrt{R_sR_d}\), \(R_x=\beta\sqrt{R_sR_d}\) and \(0<\alpha,\beta\leq 1.\) If not only \(\alpha=1,\) but also \(\beta=1,\) it is already known from (\ref{eq: optimalEqualArm}) that \(P(1,R_2,R_x)\equiv 1\) independently of \(r.\)

\section{Response of the optimal Wheatstone bridge to resistance variations}
Once the Wheatstone bridge has been balanced, in many practical scenarios, the adjusted bridge configuration is kept fixed and serves as the reference state \cite{Pallas2001}. This is particularly relevant in sensor applications, where changes in the resistance \(R_x\) caused by the measurand are detected as variations in the bridge output voltage \(V_d\). In this case, the bridge is operated in deflection mode rather than in null-balance mode: the bridge is not readjusted after the initial balancing, and the resulting output voltage directly reflects the deviation \(\Delta\) from its reference value \(R_x\). Accordingly, in this section, the optimal configuration derived in the preceding sections is considered under deflection operation, and the corresponding explicit functional relationship is derived below.

To this end, let \(\Delta\) denote a perturbation of \(R_x\) resulting in an absolute change in resistance \(R_x+\Delta.\) Note that alternatively, a relative change in \(R_x\) can analogously be considered by replacing \(\Delta\) with \(\Delta R_x.\)  According to (\ref{eq: R}), the corresponding matrix representation of the system is given by \(\bm{R_\Delta} \bm{I_\Delta}=\bm{V}\) with \(\bm{R_\Delta}=\bm{R}+\bm{\Delta},\) where
\[\bm{\Delta} = \begin{pmatrix}
0 & 0 & 0 \\
0 & 0 & 0 \\
 - \Delta & \Delta & - \Delta
\end{pmatrix}.\]
Furthermore, let \(\bm{R_{\Delta,V}}\) denote the corresponding Cramer matrix formed by replacing the last column of \(\bm{R_\Delta}\) by
\(\bm{V}=\left(
  \begin{array}{c}
    V_s,~0,~0
  \end{array}
\right)^T,\)
i.e., 
\[\bm{R_{\Delta,V}} = \begin{pmatrix}
R_{s} & R_{1} + R_{2} & V_{s} \\
 - R_{3} & R_{1} + R_{3} & 0 \\
 - {(R}_{x} + \Delta) & R_{x} + \Delta + R_{2} & 0
\end{pmatrix}.\]
By using (\ref{eq: CramerDet}), it follows
\[\det \left(\bm{R_{\Delta,V}}\right) = \det\left(\bm{R_{V}}\right) + V_{s}R_{1}\Delta,\]
which implies 
\[\det\left(\overline{\bm{R_{\Delta,V}}}\right) =V_{s}kR_{2}\Delta\]
for the balanced state \(R_1R_x=R_2R_3,\) which was subsequently perturbed by \(\Delta.\)
Let \(V_{d,\Delta}\) denote the corresponding bridge output voltage across \(R_d\) for the balanced system perturbed by \(\Delta.\) Leveraging Cramer's rule analogous to (\ref{eq: Id_Cramer}), it can be deduced that
\begin{align*}
  V_{d,\Delta}=\frac{R_{d}\det\left( \overline{\bm{R_{\Delta,V}}} \right)}{\det\left( \overline{\bm{R_{\Delta}}} \right)} = \frac{V_{s}R_{d}kR_{2}\Delta}{\det\left( \overline{\bm{R_{\Delta}}} \right)}=\frac{V_{s}R_{d}kR_{2}\Delta}{\det\left( \bm{\overline{R}}+\bm{\Delta} \right)}.
\end{align*}

The rank-one matrix \(\bm{\Delta}\) can be written as
\begin{align*}
  \bm{\Delta}=\Delta\begin{pmatrix}
0 \\
0 \\
1
  \end{pmatrix}
  \begin{pmatrix}
-1 & 1 & -1
\end{pmatrix},
\end{align*}
such that the matrix determinant lemma yields
\begin{align*}
  \det\left( \overline{\bm{R}}+\bm{\Delta}\right)=\det\left( \overline{\bm{R}}\right)+\Delta  \begin{pmatrix}
    -1 & 1 & -1
  \end{pmatrix}
  \text{adj}(\overline{\bm{R}})\begin{pmatrix}
0 \\
0 \\
1
  \end{pmatrix},
\end{align*}
where \(\text{adj}(\overline{\bm{R}})\) denotes the adjugate of \(\overline{\bm{R}}\), cf. \cite{Harville2008}.

Therefore, it holds true that
\begin{align*}
  V_{d,\Delta}=\frac{V_{s}R_{d}kR_{2}\Delta}{\det\left( \overline{\bm{R}}\right)+\Delta  \begin{pmatrix}
    -1 & 1 & -1
  \end{pmatrix}
  \text{adj}(\overline{\bm{R}})\begin{pmatrix}
0 \\
0 \\
1
  \end{pmatrix}}
\end{align*}
and after rearranging ultimately
\begin{equation}\label{eq: DeltaExplicit}
  \Delta(V_{d,\Delta}) = \frac{\det\left( \overline{\bm{R}} \right)}{\frac{V_{s}}{V_{d,\Delta}}{R_{d}kR}_{2} + \begin{pmatrix}
  1 & -1 & 1
\end{pmatrix}\text{adj}\left( \overline{\bm{R}} \right)\begin{pmatrix}
0 \\
0 \\
1
\end{pmatrix}}.
\end{equation}
By definition of the adjugate, it follows
\begin{align*}
  \text{adj}\left( \overline{\bm{R}} \right) \cdot \begin{pmatrix}
0 \\
0 \\
1
\end{pmatrix} =
\begin{pmatrix}
\det\begin{pmatrix}\overline{R}_{12} & \overline{R}_{13}\\ \overline{R}_{22} & \overline{R}_{23}\end{pmatrix} \\
-\det\begin{pmatrix}\overline{R}_{11} & \overline{R}_{13}\\ \overline{R}_{21} & \overline{R}_{23}\end{pmatrix} \\
\det\begin{pmatrix}\overline{R}_{11} & \overline{R}_{12}\\ \overline{R}_{21} & \overline{R}_{22}\end{pmatrix}
\end{pmatrix},
\end{align*}
which can explicitly be evaluated as
\begin{align*}
   \text{adj}\left( \overline{\bm{R}} \right) \cdot \begin{pmatrix}
0 \\
0 \\
1
\end{pmatrix}=
  \begin{pmatrix}
R_{2}\left( (1 + k)R_{d} + k\left( R_{2} + R_{x} \right) \right) \\
 - R_{s}R_{d} + kR_{2}R_{x} \\
k\left( R_{s}\left( R_{2} + R_{x} \right) + (1 + k)R_{2}R_{x} \right)
\end{pmatrix}
\end{align*}
by using the representation of \(\overline{\bm{R}}\) given in (\ref{eq: Rbalanced}). Hence, it can be deduced that
\begin{equation}\label{eq: adj_R}
  \begin{split}
  &\begin{pmatrix}
  1 & -1 & 1
\end{pmatrix} \cdot \text{adj}\left( \overline{\bm{R}} \right) \cdot \begin{pmatrix}
0 \\
0 \\
1
\end{pmatrix}\\
  &\quad=R_{2}\left( (1 + k)R_{d} + k\left( R_{2} + R_{x} \right) \right) + R_{s}R_{d}\\
  &\qquad - kR_{2}R_{x}+ k\left( R_{s}\left( R_{2} + R_{x} \right) + (1 + k)R_{2}R_{x} \right)\\
  &\quad = (1 + k)R_{2}R_{d} + kR_{2}^{2} + R_{s}R_{d} \\
  &\qquad+ k\left( R_{s}\left( R_{2} + R_{x} \right) + (1 + k)R_{2}R_{x} \right)
 \end{split}
\end{equation}

Therefore, (\ref{eq: DeltaExplicit}) expresses the resistance change \(\Delta\) as a function of the corresponding output voltage change \(V_{d,\Delta}\) for an initially balanced Wheatstone bridge satisfying \(R_1R_x=R_2R_3,\) where \(k=\frac{R_1}{R_2}.\) Furthermore, the involved determinant and adjugate terms can explicitly be evaluated with the help of (\ref{detBalancedR}) and (\ref{eq: adj_R}), respectively.

The functional relationship for \(\Delta\) presented in (\ref{eq: DeltaExplicit}) is valid for any Wheatstone bridge configuration determined by the choice of \(k\) and \(R_2.\) In the following, the explicit functional relationship is derived for the optimal Wheatstone bridge configuration characterized by \(k^{opt}\) and \(R_2^{opt}\).

\begin{figure*}[!t]
\normalsize
\begin{equation}\label{eq: adjOpt}
  \begin{split}
      &\begin{pmatrix}
1 & - 1 & 1
\end{pmatrix} \cdot \text{adj}\left(\overline{\bm{R}}^{opt}\right) \cdot \begin{pmatrix}
0 \\
0 \\
1
\end{pmatrix} \\
  &\quad= R_{2}^{opt}R_{d} + \sqrt{R_{s}R_{d}}R_{d} + \sqrt{R_{s}R_{d}}R_{2}^{opt} + R_{s}R_{d}+ \sqrt{R_{s}R_{d}}R_{s} + k^{opt}R_{s}R_{x} + \sqrt{R_{s}R_{d}}R_{x} + k^{opt}\sqrt{R_{s}R_{d}}R_{x} \\
  &\quad= R_{s}R_{d} + R_{2}^{opt}\left( R_{d} + \sqrt{R_{s}R_{d}} \right)+ \sqrt{R_{s}R_{d}}\left( R_{s} + R_{d} + R_{x} \right) + k^{opt}R_{x}\left( R_{s} + \sqrt{R_{s}R_{d}} \right)\\
  &\quad = R_{s}R_{d} + \sqrt{\frac{R_{s}R_{x}\left( R_{d} + R_{x} \right)}{R_{s} + R_{x}}}\left( R_{d} + \sqrt{R_{s}R_{d}} \right)+ \sqrt{R_{s}R_{d}}\left( R_{s} + R_{d} + R_{x} \right)+ \sqrt{\frac{R_{d}\left( R_{s} + R_{x} \right)}{R_{x}\left( R_{d} + R_{x} \right)}}R_{x}\left( R_{s} + \sqrt{R_{s}R_{d}} \right) \\
      &\quad= R_{s}R_{d} + \sqrt{R_{s}R_{d}}\left( R_{s} + R_{d} + R_{x} \right)+ \left( R_{s}\sqrt{R_{d}} + R_{d}\sqrt{R_{s}} \right)\sqrt{R_{x}}\left( \sqrt{\frac{R_{d} + R_{x}}{R_{s} + R_{x}}} + \sqrt{\frac{R_{s} + R_{x}}{R_{d} + R_{x}}} \right) \\
      &\quad= \sqrt{R_{s}R_{d}}\left( \sqrt{R_{s}R_{d}} + R_{s} + R_{d} + R_{x} + \frac{\left( \sqrt{R_{s}} + \sqrt{R_{d}} \right)\sqrt{R_{x}}\left( R_{s} + {2R}_{x} + R_{d} \right)}{\sqrt{\left( R_{s} + R_{x} \right)\left( R_{d} + R_{x} \right)}} \right)
  \end{split}
\end{equation}
\hrulefill 
\vspace*{4pt} 
\end{figure*}

Using the representation of \(\det\left(\overline{\bm{R}}\right)\) and the explicit evaluation of the adjugate term for the optimal setup given in (\ref{eq: detOpt}) and (\ref{eq: adjOpt}), respectively, the functional relationship (\ref{eq: DeltaExplicit}) for the optimal setup becomes
\begin{align*}
  &\Delta(V_{d,\Delta})= - \frac{R_{\Sigma}^{2}(R_s,R_d,R_x)}{\frac{V_{s}}{V_{d,\Delta}}R_{d} + R_{\Gamma}(R_s,R_d,R_x)}
\end{align*}
with
\begin{align*}
  &R_\Sigma(R_s,R_d,R_x)\\
  &\quad=\sqrt{R_{s}R_{x}} + \sqrt{\left( R_{s} + R_{x} \right)\left( R_{d} + R_{x} \right)} + \sqrt{{R_{d}R}_{x}}
\end{align*}
and
\begin{align*}
  R_{\Gamma}(R_s,R_d,R_x)&=\sqrt{R_{s}R_{d}} + R_{s} + R_{d} + R_{x} \\
                         &+ \frac{\left( \sqrt{R_{s}} + \sqrt{R_{d}} \right)\sqrt{R_{x}}\left( R_{s} + {2R}_{x} + R_{d} \right)}{\sqrt{\left( R_{s} + R_{x} \right)\left( R_{d} + R_{x} \right)}}.
\end{align*}

\section{Conclusion}
This paper provides the first analytical solution for the optimal resistance configuration for maximizing null sensitivity of the Wheatstone bridge with finite source and detector resistances. From a scientific perspective, this result resolves a longstanding open challenge in the academic literature, outperforming the traditional equal-arm setup and eliminating the need for numerical approximations.

From a practical perspective, the proposed closed-form optimal solutions (\ref{eq: R2Opt}) and (\ref{eq: R1Opt}) for the fixed bridge components offer an immediate impact on high-precision measurement applications, enabling engineers to optimize sensitivity directly through the bridge design. As demonstrated in Section~\ref{sec: comparison}, the optimal configuration outperforms the classical equal-arm design, particularly when the ratio of detector to source resistance \(\frac{R_d}{R_s}\) is not sufficiently large, or when \(R_x\) deviates significantly from the equal-arm sweet spot \(\sqrt{R_sR_d}\) and is instead comparable in magnitude to \(R_s\) or \(R_d.\)


\vfill

\end{document}